\documentclass[10pt,draft,reqno]{amsart}
     \makeatletter
     \def\section{\@startsection{section}{1}%
     \z@{.7\linespacing\@plus\linespacing}{.5\linespacing}%
     {\bfseries
     \centering
     }}
     \def\@secnumfont{\bfseries}
     \makeatother
\newtheorem{theorem}{Theorem}[section]
\newtheorem{lemma}[theorem]{Lemma}
\newtheorem{proposition}[theorem]{Proposition}

\theoremstyle{definition}
\newtheorem{definition}[theorem]{Definition}

\theoremstyle{remark}
\newtheorem{remark}[theorem]{Remark}
\numberwithin{equation}{section} 

\newcommand{\R}{\mathbb{R}}
\newcommand{\Tr}{\mathop{\mathrm{Tr}}}
\newcommand{\Div}{\mathop{\mathrm{Div}}}
\renewcommand{\d}{\/\mathrm{d}\/}

\begin{document}
\title[Stochastic Navier-Stokes Equation with Artificial Compressibility]
{Stochastic $2$-D Navier-Stokes Equation with Artificial
Compressibility}

\author[Utpal Manna]{Utpal Manna*}
\thanks{* This research is supported by Army Research Office, Probability and
Statistics Program, grant number DODARMY1736}
\address{Utpal Manna: Department of Mathematics, University of Wyoming, Laramie, WY 82071, USA}
\email{utpal@uwyo.edu} \urladdr{http://uwyo.edu/utpal/}

\author{J.L. Menaldi}
\address{J.L. Menaldi: Department of Mathematics, Wayne State University, Detroit, MI 48202, USA}
\email{jlm@math.wayne.edu}
\urladdr{http://www.math.wayne.edu/~menaldi/}

\author[S.S. Sritharan]{S.S. Sritharan*}
\address{S.S. Sritharan: Department of Mathematics, University of Wyoming, Laramie, WY 82071, USA}
\email{sri@uwyo.edu}
\urladdr{http://math.uwyo.edu/Prof\_home/Sri.html}



\begin{abstract}
In this paper we study the stochastic Navier-Stokes equation with
artificial compressibility. The main results of this work are the
existence and uniqueness theorem for strong solutions and the limit
to incompressible flow. These results are obtained by utilizing a
local monotonicity property of the sum of the Stokes operator and
the nonlinearity.
\end{abstract}

\maketitle

\section{Introduction}

The stochastic Navier-Stokes equation is a well accepted model for
atmospheric, aero and ocean dynamics. Chandrasekhar~\cite{Cs89} and
Novikov~\cite{Nv65} first studied the Navier-Stokes equation with
external random forces. After that several approaches have been
poposed, from the classic paper by Bensoussan and
Temam~\cite{BeTe73} to some more recent results, e.g., by
Bensoussan~\cite{Be95}, by Flandoli and Gatarek~\cite{FfDg95} and by
Menaldi and Sritharan~\cite{MeSr00}.

This paper is concerned with the existence and uniqueness of strong
solutions for the Stochastic 2-D Navier-Stokes equation with
artificial compressibility in bounded domains. The concept of
artificial compressibility was first introduced by
Chorin~\cite{Ch68,Ch69} and Temam~\cite{Te691,Te692}, in order to
overcome the computational difficulties connected with the
incompressibility constraint. Using the classical Sobolev
compactness embedding and exploiting the classical Lions~\cite{Li59}
method of fractional derivatives, Temam in his
papers~\cite{Te691,Te692} and in his book~\cite{Te95}(chapter 3)
proved the existence, uniqueness and convergence of the
deterministic Navier-Stokes equation with artificial compressibility
in bounded domains.

In the rest of this section we formulate the abstract Navier-Stokes
problem with artificial compressibility. We describe some standard
well known results including the local monotonicity property of the
Navier-Stokes operator. In Section 2 we establish certain new a
priori estimates involving exponential weight for stochastic
Navier-Stokes equation with artificial compressibility. These
estimates play a fundamental role in the proof of existence and
uniqueness of strong solutions proved in the second half of Section
2. The monotonicity argument used here is the generalization of the
classical Minty-Browder method for dealing with local monotonicity.
This method was first used by Menaldi and Sritharan~\cite{MeSr00}
and for multiplicative noise by Sritharan and Sundar~\cite{SrSu06}.
For similar ideas see also Barbu and Sritharan~\cite{BaSr98,
BaSr01}. In the last part of Section 2 we discuss the convergence of
the corresponding perturbed problem. Here the use of local
monotonicity avoids the classical method based on compactness and
thus the results apply to unbounded domains and hence the existence
and the uniqueness as well as convergence to incompressible flow are
new even in the deterministic case.


\section{Abstract Mathematical Framework and Local Monotonicity}

Let $\mathcal{O}\subset\mathbb{R}^{2}$ be a bounded domain (for the
sake of simplicity) with smooth boundary, $\mathbf{u}$ the velocity
and $p$ the pressure fields. The Navier-Stokes problem (with
Newtonian constitutive relationship and \emph{artificial
compressible} medium) can be written as follows
\begin{equation}\label{e1.1}\left\{\begin{aligned}
  &\partial _{t}\mathbf{u}-\nu\triangle\mathbf{u}+(\mathbf{u}
  \cdot\nabla)\mathbf{u}+\frac{1}{2}(\Div\mathbf{u})\mathbf{u}
  +\nabla p =\mathbf{f}\; \text{ in }\;
  \mathrm{L}^2(0,T;\mathbb{H}^{-1}(\mathcal{O})),\\
  &\varepsilon\partial_{t}p + \Div\mathbf{u}=0 \;
  \text{ in }\; \mathrm{L}^2(0,T;\mathrm{L}^2(\mathcal{O})),
\end{aligned}\right.\end{equation}
with the initial conditions
\begin{equation}
  \mathbf{u}(0)=\mathbf{u}_{0}\quad
  \text{in}\quad\mathbb{L}^2(\mathcal{O})\qquad\text{and}\qquad
  p(0)=p_0\quad\text{in}\quad \mathrm{L}^2(\mathcal{O}),
\end{equation}
where $\varepsilon>0$ is a vanishing parameter, $\mathbf{u}_{0}$
belong to
$\mathbb{L}^2(\mathcal{O})=\mathrm{L}^2(\mathcal{O},\R^2),$ $\nu$ is
the kinematic viscosity, $p$ denotes pressure and is a scalar-valued
function and the (force) field $\mathbf{f}$ is in
$\mathrm{L}^2(0,T;\mathbb{L}^2(\mathcal{O})).$  A solution
$(\mathbf{u},p)$ should belongs to the space
$\mathrm{L}^2(0,T;\mathbb{H}^1_0(\mathcal{O})\times
\mathrm{L}^2(\mathcal{O})),$ with
$\mathbb{H}^1_0(\mathcal{O})=\mathrm{H}^1_0(\mathcal{O},\R^2)$ and
$\mathbb{H}^{-1}(\mathcal{O})$ its dual space.  The second equation
in \eqref{e1.1} is an artificial state equation of a slightly
compressible medium and the extra term
$(1/2)(\Div\mathbf{u})\mathbf{u}$ is a stabilization term to handle
the nonlinearity. The standard spaces used here are as follows: \\
$\mathbb{H}^1_0(\mathcal{O})$ with the norm
\begin{equation}
\|\mathbf{v}\|_{_{\mathbb{H}^1_0}} := \Big(\int_\mathcal{O}
|\nabla\mathbf{v}|^2 dx \Big)^{1/2} = \|\mathbf{v}\|,
\end{equation}
and $\mathbb{L}^2(\mathcal{O})$ with the norm
\begin{equation}
\|\mathbf{v}\|_{_{\mathbb{L}^2}} := \Big(\int_\mathcal{O}
|\mathbf{v}|^2 dx \Big)^{1/2} = |\mathbf{v}|.
\end{equation}

 Using the Gelfand triple $\mathbb{H}^1_0(\mathcal{O}) \subset
\mathbb{L}^2(\mathcal{O})\subset \mathbb{H}^{-1}(\mathcal{O})$ we
may consider $\triangle$  or $\nabla$ as a linear map from
$\mathbb{H}^1_0(\mathcal{O})$ or $\mathrm{L}^2(\mathcal{O})$ into
the dual of $\mathbb{H}^1_0(\mathcal{O})$ respectively. The inner
product in the $\mathbb{L}^{2}$ or $\mathrm{L}^{2}$ is denoted by
$(\cdot,\cdot)$ and the induced duality by
$\langle\cdot,\cdot\rangle.$  Thus, for any $\mathbf{u} = (u_i),$
$\mathbf{v} = (v_i)$ and $\mathbf{w} = (w_i)$ in
$\mathbb{H}^1_0(\mathcal{O})$ and $p$ in $\mathrm{L}^2(\mathcal{O})$
we have
\begin{equation}
  \langle -\nu\triangle\mathbf{u},\mathbf{w}\rangle = \nu\sum_{i,j}
  \int_\mathcal{O} \partial_i u_j \, \partial_i w_j \d x,
\end{equation}
\begin{equation}\label{e1.6}
  \langle - \nabla p,\mathbf{w}\rangle = -\sum_i
  \int_\mathcal{O} \partial_i p \, w_i \d x =
  \int_\mathcal{O}  p \,\partial_i w_i \d x =
  \langle p,\Div\mathbf{w}\rangle
\end{equation}
and
\begin{equation}
  \langle (\mathbf{u}\cdot\nabla)\mathbf{v}, \mathbf{w}\rangle = \sum_{i,j}
  \int_\mathcal{O} u_i\, \partial_i v_j\, w_j \d x.
\end{equation}
It is clear that $\mathbf{u}\mapsto\Div\mathbf{u}$ is a linear
continuous operator from $\mathbb{H}^1_0(\mathcal{O})$ into
$\mathrm{L}^2(\mathcal{O}).$  Next, an integration by parts and
H\"older inequality yields
\begin{align}\label{e1.8}
  &\langle (\mathbf{u}\cdot\nabla)\mathbf{v}, \mathbf{w}\rangle =
  - \langle (\Div\mathbf{u})\mathbf{w},\mathbf{v}\rangle
  - \langle (\mathbf{u}\cdot\nabla)\mathbf{w}, \mathbf{v}\rangle, \\
  &|\langle (\mathbf{u}\cdot\nabla)\mathbf{v}, \mathbf{w}\rangle | \leq
  C\sum_{i,j}\|u_i w_j\|_{_{\mathrm{L}^2(\mathcal{O},\R^2)}} \|\partial_i
  v_j\|_{_{\mathrm{L}^2(\mathcal{O},\R^2)}}, \label{e1.9}
\end{align}
and in the right-hand-side we can use $\mathbb{L}^4$-norms to
estimate the product $u_i v_j.$


\begin{lemma}
For any real-valued smooth functions $\varphi$ and $\psi$ with
compact support in $\R^2,$ the following hold:
\begin{align}\label{e1.10}
  & \| \varphi\,\psi \|_{_{\mathrm{L}^2}}^2 \leq  \|\varphi \,
  \partial_1\varphi \|_{_{\mathrm{L}^1}} \|\psi \, \partial_2\psi
  \|_{_{\mathrm{L}^1}}, \\
  & \| \varphi \|_{_{L^4}}^4 \leq 2 \|\varphi\|_{_{\mathrm{L}^2}}^2
  \|\nabla \varphi\|_{_{\mathbb{L}^2}}^2. \label{e1.11}
\end{align}
\end{lemma}

\begin{proof}
The results stated above are classical and well known~\cite{La69}.
\end{proof}

As in Temam~\cite{Te95}, chapter 3, the non-linear term is a
trilinear continuous form on
$\mathbb{H}^1_0(\mathcal{O})\times\mathbb{H}^1_0(\mathcal{O})\times\mathbb{H}^1_0(\mathcal{O})$
\begin{equation}\label{e1.12}
  \hat{b}(\mathbf{u},\mathbf{v},\mathbf{w}):=\langle\hat{B}
  (\mathbf{u},\mathbf{v}),\mathbf{w}\rangle:=\frac{1}{2}\sum_{i,j}
  \int_\mathcal{O} [u_i\, \partial_i v_j\, w_j -
  u_i\, \partial_i w_j\, v_j]\d x,
\end{equation}
where
\begin{equation}
  \hat{B}(\mathbf{u})=\hat{B}(\mathbf{u},\mathbf{u})=[(\mathbf{u}\cdot\nabla)
  + \frac{1}{2}\Div\mathbf{u}]\mathbf{u}.
\end{equation}
We have the following lemmas. \\


\begin{lemma}
Let $\mathbf{u}$ and $\mathbf{w}$ be in the spaces
$\mathbb{H}^1_0(\mathcal{O}, \mathbb{R}^{2})$ and
$\mathbb{L}^4(\mathcal{O}, \mathbb{R}^{2})$ respectively. Then the
following estimate holds:
\begin{equation}\label{e1.14}
|\langle\hat{B}(\mathbf{u}), \mathbf{w}\rangle| \leq 2
\|\mathbf{u}\|^{3/2} \ |\mathbf{u}|^{1/2} \
\|\mathbf{w}\|_{\mathbb{L}^4(\mathcal{O}, \mathbb{R}^{2})}.
\end{equation}
\end{lemma}

\begin{proof}
We observe that
\begin{equation}
\hat{b}(\mathbf{u}, \mathbf{v}, \mathbf{w}) = \sum_{i, j
=1}^{2}\int_{\mathcal{O}}u_{i}(D_{i}v_{j})w_{j}\ dx \ + \
\frac{1}{2}\sum_{j =1}^{2}\int_{\mathcal{O}} (div\
\mathbf{u})v_{j}w_{j} \ dx. \nonumber
\end{equation}
Using the H\"older inequality,
\begin{eqnarray}
|\hat{b}(\mathbf{u}, \mathbf{v}, \mathbf{w})| &\leq& \sum_{i, j
=1}^{2} \|u_{i}\|_{\mathbb{L}^4(\mathcal{O}, \mathbb{R}^{2})} \
|D_{i}v_{j}|_{\mathbb{L}^2(\mathcal{O}, \mathbb{R}^{2})} \
\|w_{j}\|_{\mathbb{L}^4(\mathcal{O}, \mathbb{R}^{2})}\nonumber \\
& + & \frac{1}{2} \sum_{j =1}^{2}\ |div \
\mathbf{u}|_{\mathbb{L}^2(\mathcal{O}, \mathbb{R}^{2})} \
\|v_{j}\|_{\mathbb{L}^4(\mathcal{O}, \mathbb{R}^{2})} \
\|w_{j}\|_{\mathbb{L}^4(\mathcal{O}, \mathbb{R}^{2})}.\nonumber
\end{eqnarray}
In Frobenius norm divergence can be estimated by gradient. Hence
\begin{equation}
|\hat{b}(\mathbf{u}, \mathbf{u}, \mathbf{w})| \leq\
\frac{3}{2}\|\mathbf{u}\|_{\mathbb{L}^4(\mathcal{O},
\mathbb{R}^{2})} \ \|\mathbf{u}\| \
\|\mathbf{w}\|_{\mathbb{L}^4(\mathcal{O}, \mathbb{R}^{2})}.\nonumber
\end{equation}
Using the equation \eqref{e1.11} in Lemma 2.1 we get
\begin{equation}
|\langle\hat{B}(\mathbf{u}), \mathbf{w}\rangle| \ = \
|\hat{b}(\mathbf{u}, \mathbf{u}, \mathbf{w})|\ \leq\ 2
\|\mathbf{u}\|^{3/2} \ |\mathbf{u}|^{1/2} \
\|\mathbf{w}\|_{\mathbb{L}^4(\mathcal{O}, \mathbb{R}^{2})}.\nonumber
\end{equation}
\end{proof}


\begin{lemma}
Let $\mathbf{u}$ and $\mathbf{v}$ be in the space
$\mathbb{H}^1_0(\mathcal{O}, \mathbb{R}^{2})$. Then the following
estimate holds:
\begin{equation}\label{e1.15}
|\langle\hat{B}(\mathbf{u})-\hat{B}(\mathbf{v}),\mathbf{u}-\mathbf{v}
  \rangle | \leq \frac{\nu}{2} \,
  \|\mathbf{u}- \mathbf{v}\|^2 + \frac{27}{2\nu^3} \,
  |\mathbf{u}-\mathbf{v}|^2 \,
  \|\mathbf{v}\|_{_{\mathbb{L}^4(\mathcal{O},\R^2)}}^4.
\end{equation}
\end{lemma}

\begin{proof}
 For any given $\mathbf{u}$ and $\mathbf{v}$ in $\mathbb{H}^1_0(\mathcal{O}, \mathbb{R}^{2})$
 , we have from \eqref{e1.12}
\begin{align}\label{1.16}
\langle \hat{B}(\mathbf{u}), \mathbf{u}\rangle = 0,
\end{align}
and
\begin{align}\label{1.17}
\langle \hat{B}(\mathbf{u},\mathbf{v}), \mathbf{v}\rangle
=\langle(\mathbf{u}\cdot\nabla)\mathbf{v}, \mathbf{v}\rangle +
\frac{1}{2}\langle (\Div\mathbf{u})\mathbf{v},\mathbf{v}\rangle = 0.
\end{align}
Then using \eqref{1.16} we obtain
\begin{align}
 \langle \hat{B}(\mathbf{u})-\hat{B}(\mathbf{v}),
  \mathbf{u}-\mathbf{v}\rangle &= \langle \hat{B}(\mathbf{u}), \mathbf{u}\rangle -\langle
  \hat{B}(\mathbf{u}),\mathbf{v}\rangle + \langle \hat{B}(\mathbf{v}), \mathbf{v}\rangle
  - \langle \hat{B}(\mathbf{v}),\mathbf{u}\rangle\nonumber\\
  &= -\langle\hat{B}(\mathbf{u}),\mathbf{v}\rangle - \langle
\hat{B}(\mathbf{v}),\mathbf{u}\rangle.\label{1.18}
\end{align}
Now
\begin{align}
\langle \hat{B}(\mathbf{u}-\mathbf{v}), \mathbf{v}\rangle &=
\langle((\mathbf{u}-\mathbf{v})\cdot\nabla)(\mathbf{u}-\mathbf{v}) +
\frac{1}{2}(\Div(\mathbf{u}-\mathbf{v}))(\mathbf{u}-\mathbf{v}),
\mathbf{v}\rangle\nonumber\\
&=\langle\hat{B}(\mathbf{u}),\mathbf{v}\rangle +
\langle\hat{B}(\mathbf{v}),\mathbf{v}\rangle - \langle
(\mathbf{u}\cdot\nabla)\mathbf{v} +
\frac{1}{2}(\Div\mathbf{u})\mathbf{v},\mathbf{v}\rangle\nonumber\\
&\quad - \langle (\mathbf{v}\cdot\nabla)\mathbf{u} +
\frac{1}{2}(\Div\mathbf{v})\mathbf{u},\mathbf{v}\rangle.\nonumber
\end{align}
Hence using \eqref{1.16} and \eqref{1.17} we have
\begin{align}
\langle \hat{B}(\mathbf{u}-\mathbf{v}), \mathbf{v}\rangle
&=\langle\hat{B}(\mathbf{u}),\mathbf{v}\rangle - \langle
(\mathbf{v}\cdot\nabla)\mathbf{u} +
\frac{1}{2}(\Div\mathbf{v})\mathbf{u},\mathbf{v}\rangle\nonumber\\
&=\langle\hat{B}(\mathbf{u}),\mathbf{v}\rangle  +
\langle\hat{B}(\mathbf{v}),\mathbf{u}\rangle - \langle
(\mathbf{v}\cdot\nabla)\mathbf{u} +
\frac{1}{2}(\Div\mathbf{v})\mathbf{u},\mathbf{v}\rangle\nonumber\\
&\quad - \langle (\mathbf{v}\cdot\nabla)\mathbf{v} +
\frac{1}{2}(\Div\mathbf{v})\mathbf{v},\mathbf{u}\rangle\nonumber\\
&=\langle\hat{B}(\mathbf{u}),\mathbf{v}\rangle  +
\langle\hat{B}(\mathbf{v}),\mathbf{u}\rangle \nonumber\\
&\quad - \langle (\mathbf{v}\cdot\nabla)(\mathbf{u}+\mathbf{v}) +
\frac{1}{2}(\Div\mathbf{v})(\mathbf{u}+\mathbf{v}),\mathbf{u}+\mathbf{v}\rangle\nonumber\\
&\quad + \langle (\mathbf{v}\cdot\nabla)\mathbf{u} +
\frac{1}{2}(\Div\mathbf{v})\mathbf{u},\mathbf{u}\rangle\nonumber\\
&\quad + \langle (\mathbf{v}\cdot\nabla)\mathbf{v} +
\frac{1}{2}(\Div\mathbf{v})\mathbf{v},\mathbf{v}\rangle.\nonumber
\end{align}
With the help of \eqref{1.16} and \eqref{1.17} last three terms of
the right hand side vanish.Thus
\begin{align}\label{1.19}
\langle \hat{B}(\mathbf{u}-\mathbf{v}), \mathbf{v}\rangle =
\langle\hat{B}(\mathbf{u}),\mathbf{v}\rangle  +
\langle\hat{B}(\mathbf{v}),\mathbf{u}\rangle.
\end{align}
 Thus \eqref{1.18} and \eqref{1.19} yield
\begin{align}
\langle \hat{B}(\mathbf{u})-\hat{B}(\mathbf{v}),
  \mathbf{u}-\mathbf{v}\rangle = -\langle \hat{B}(\mathbf{u}-\mathbf{v}),
  \mathbf{v}\rangle,
\end{align}
which by Lemma 2.2 gives the estimate
\[
  | \langle \hat{B}(\mathbf{u}) - \hat{B}(\mathbf{v}), \mathbf{u}-
  \mathbf{v}\rangle | \leq 2\, \|\mathbf{u}-\mathbf{v}\|^{3/2} \,
  |\mathbf{u}-\mathbf{v}|^{1/2} \,
  \|\mathbf{v}\|_{_{\mathbb{L}^4(\mathcal{O},\R^2)}},
\]
where $\|\cdot\|$ and $|\cdot|$ denotes the norm in
$\mathbb{H}^1_0(\mathcal{O})$ and $\mathbb{L}^2(\mathcal{O})$
respectively. \\
Now using the fact that for any two real numbers $a$, $b$ and any
$p$, $q$ $> 1$ with $\frac{1}{p} + \frac{1}{q} = 1$,
$$ ab \leq \frac{|a|^{p}}{p} + \frac{|b|^{q}}{q} ,$$
we obtain the estimate \eqref{e1.15}.
\end{proof}

Notice that $\mathbb{H}^1_0(\mathcal{O})$ is continuously included
in $\mathbb{L}^4(\mathcal{O})=\mathrm{L}^4(\mathcal{O},\R^2)$ and
$\mathbf{u}\mapsto(\nabla\cdot\mathbf{u})\mathbf{u}$ is a
(nonlinear) continuous mapping from $\mathbb{H}^1_0(\mathcal{O})$
into its dual $\mathbb{H}^{-1}(\mathcal{O}).$ Hence the nonlinear
operator $\hat{B}(\cdot)$ can be considered as a map from
$\mathbb{H}^1_0(\mathcal{O})$ into the space
$\mathbb{H}^{-1}(\mathcal{O})\cap \mathbb{L}^{4/3}(\mathcal{O},
\mathbb{R}^2)$. Then combination of previous lemmas yield the
following local monotonicity property.


\begin{lemma}
For a given $r > 0$, let us denote by $\mathbb{B}_{r}$ the (closed)
$\mathbb{L}^{4}$-ball in $\mathbb{H}^1_0$
\begin{equation}
  \mathbb{B}_{r}=\{\mathbf{v}\in \mathrm{H}^1_0(\mathcal{O},\mathbb{R}^{2})
  \text{ ; }\Vert \mathbf{v}
  \Vert _{_{\mathrm{L}^{4}(\mathcal{O},\mathbb{R}^{2})}}\leq r\},
\end{equation}
then the nonlinear operator $\mathbf{u}\mapsto
A\mathbf{u}+\hat{B}(\mathbf{u}):=
-\nu\triangle\mathbf{u}+[(\mathbf{u}\cdot\nabla)+(1/2)\Div\mathbf{u}]\mathbf{u}$
is monotone in the convex ball $\mathbb{B}_{r}$ i.e.,
\begin{equation}\label{e1.17}
  \langle A\mathbf{w},\mathbf{w}\rangle +\langle \hat{B}(\mathbf{u})
  -\hat{B}(\mathbf{v}),  \mathbf{w}\rangle +\frac{27r^{4}}{2\nu ^{3}}
  \,|\mathbf{w}|^{2}\geq \frac{\nu}{2}\,\Vert \mathbf{w}\Vert^{2},
\end{equation}
$\forall \mathbf{u}\in \mathrm{H}^1_0(\mathcal{O},\mathbb{R}^{2}),$
$\mathbf{v}\in \mathbb{B}_{r}$ and
$\mathbf{w}=\mathbf{u}-\mathbf{v}.$
\end{lemma}

\begin{proof}
First, it is clear that
$$\langle A\mathbf{w},\mathbf{w}\rangle = \nu \|\mathbf{w}\|^{2},$$
and the equation\eqref{e1.15} yields
$$\langle \hat{B}(\mathbf{u})-\hat{B}(\mathbf{v}),
\mathbf{w}\rangle \geq - \frac{\nu}{2}\,\Vert \mathbf{w}\Vert^{2} -
\frac{27r^{4}}{2\nu ^{3}} \,|\mathbf{w}|^{2}.$$ Summing these
equations up we get the desired result \eqref{e1.17}.
\end{proof}


\section{Stochastic $2$-D Navier-Stokes Equation with \qquad\\ Artificial Compressibility}

Let us consider the Navier-Stokes equation subject to a random
(Gaussian) term i.e., the forcing field $\mathbf{f}$ has a mean
value still denoted by $\mathbf{f}$ and a noise denoted by
$\dot{\mathbf{G}}.$ We can write (to simplify notation we use
time-invariant forces) $\mathbf{f}(t)=\mathbf{f}(x,t)$ and the noise
process $\dot{\mathbf{G}}(t)=\dot{\mathbf{G}}(x,t)$ as a series
$d\mathbf{G}_{k}=\sum_{k}\mathbf{g}_{k}(x,t)dw_{k}(t),$ where
$\mathbf{g}=(\mathbf{g}_{1},\mathbf{g}_{2},\cdots )$ and
$w=(w_{1},w_{2},\ldots )$ are regarded as $\ell^{2}$-valued
functions in $x$ and $t$ respectively. The stochastic noise process
represented by
$\mathbf{g}(t)dw(t)=\sum_{k}\mathbf{g}_{k}(x,t)dw_{k}(t,\omega )$ is
normal distributed in $\mathbb{H}$ with a trace-class co-variance
operator denoted by $\mathbf{g}^{2}=\mathbf{g}^{2}(t)$ and given by
\begin{equation}
  \left\{\begin{split}
  & (\mathbf{g}^{2}(t)\mathbf{u},\mathbf{v})=\sum_{k}(\mathbf{g}_{k}(t),%
  \mathbf{u})\,(\mathbf{g}_{k}(t),\mathbf{v}), \\
  & \mathrm{Tr}(\mathbf{g}^{2}(t))=\sum_{k}|\mathbf{g}_{k}(t)|^{2}<\infty .
  \end{split}\right.
\end{equation}

We interpret the stochastic Navier-Stokes equations as an It\^o
stochastic equations in variational form
\begin{equation}
  \left\{\begin{split}\label{e2.2}
  &\d(\mathbf{u}(t),\mathbf{v}) + \langle -\nu\triangle\mathbf{u}(t) +
  [(\mathbf{u}(t)\cdot\nabla) +
  \frac{1}{2}\Div\mathbf{u}(t)]\mathbf{u}(t)\\
  &\qquad\qquad\qquad\qquad +\nabla p(t), \mathbf{v} \rangle \, \d t
  = ( \mathbf{f},\mathbf{v} ) \, \d t + \sum_k (
  \mathbf{g}_k,\mathbf{v})\, \d w_k(t),\!\!\\
  &\langle \varepsilon \dot{p}(t)+ \Div\mathbf{u}(t),q\rangle = 0,
  \end{split}\right.
\end{equation}
in $(0,T),$ with the initial condition
\begin{equation}
  (\mathbf{u}(0),\mathbf{v}) =
  (\mathbf{u}_0,\mathbf{v})\quad\text{and}\quad
  (p(0),q)=(p_0,q),
\end{equation}
for any $\mathbf{v}$ in the space $\mathbb{H}^1_0(\mathcal{O})$ and
any $q$ in $\mathrm{L}^2(\mathcal{O}).$

A finite-dimensional (Galerkin) approximation of the stochastic
Navier-Stokes equation can be defined as follows. Let $\{\mathbf{e}_{1},%
\mathbf{e}_{2},\ldots \}$ be a complete orthonormal system (i.e., a
basis) in the Hilbert space $\mathbb{L}^2(\mathcal{O})$ belonging to
the space $\mathbb{H}^1_0(\mathcal{O})$ (and
$\mathbb{L}^4(\mathcal{O})$). Denote by
$\mathbb{L}^2_{n}(\mathcal{O})$ the $n$-dimensional subspace of
$\mathbb{L}^2(\mathcal{O})$ and $\mathbb{H}^1_0(\mathcal{O})$ of all
linear combinations of the first $n$ elements
$\{\mathbf{e}_{1},\mathbf{e}_{2},\ldots ,\mathbf{e}_{n}\}.$ Also
denote by
$\mathrm{L}^2_{n}(\mathcal{O}):=\nabla\cdot\mathbb{L}^2_{n}(\mathcal{O})$
the image of $\nabla.$

Consider the following stochastic ODE in $\mathbb{R}^{n}$
\begin{equation}\left\{\begin{split}
  &\d(\mathbf{u}^n(t),\mathbf{v}) + \langle -\nu\triangle\mathbf{u}^n(t)
  + [(\mathbf{u}^n(t)\cdot\nabla) + \frac{1}{2}
  \Div\mathbf{u}^n(t)]\mathbf{u}^n(t) \\
  &\qquad\qquad\qquad\qquad + \nabla p^n(t),  \mathbf{v} \rangle  \d t
  = ( \mathbf{f},\mathbf{v} ) \d t + \sum_k (
  \mathbf{g}_k,\mathbf{v}) \d w_k(t),\\
  &\langle \varepsilon \dot{p}^n(t)+ \Div\mathbf{u}^n(t),q\rangle = 0,
  \end{split}\right. \label{e2.4}
\end{equation}
in $(0,T),$ with the initial condition
\begin{equation}
  (\mathbf{u}(0),\mathbf{v})=(\mathbf{u}_{0},\mathbf{v}),
\end{equation}
for any $\mathbf{v}$ in the space $\mathbb{L}^2_{n}(\mathcal{O})$
and $q$ in $\mathrm{L}^2_{n}(\mathcal{O}).$  The coefficients
involved are locally Lipschitz and we need some \emph{a priori}
estimate to show global existence of a solution $\mathbf{u}^{n}(t)$
as an adapted process in the space
$C^{0}(0,T,\mathbb{L}^2_{n}(\mathcal{O})).$


\begin{proposition}[energy estimate]
Under the above mathematical setting let
\begin{equation}
  \mathbf{f}\in \mathrm{L}^2(0,T;\mathbb{L}^2(\mathcal{O})),\;
  \mathbf{g}\in \mathrm{L}^2(0,T;\ell_{2}(\mathbb{L}^2(\mathcal{O}))),
  \; \mathbf{u}_{0}\in \mathbb{L}^2(\mathcal{O}),\;
  p_0\in \mathrm{L}^2(\mathcal{O}).\!\!\!
\end{equation}
Let $\mathbf{u}^{n}(t)$ be an adapted process in
$C^{0}(0,T,\mathbb{H}_{n})$ which solves the stochastic ODE
(\ref{e2.4}). Then we have the energy equality
\begin{equation}\left\{\begin{split}\label{e2.7}
  &\d\big[ |\mathbf{u}^{n}(t)|^2+\varepsilon|p^n(t)|^2\big]+2\nu \,
  |\nabla\mathbf{u}^{n}(t)|^{2}\d t \\
  &\qquad=\big[2\,(\mathbf{f}(t), \mathbf{u}^{n}(t))+
  \Tr(\mathbf{g}^{2}(t))\big]\d t+
  2\sum_{k}(\mathbf{g}_{k}(t),\mathbf{u}^{n}(t)) \d w_{k}(t),
\end{split}\right.\end{equation}
which yields the following estimate for any $\delta >0$
\begin{equation}\left\{\begin{split}\label{e2.8}
  E\big\{|\mathbf{u}^{n}(t)|^2 &+
  \varepsilon |p^{n}(t)|^2\big\} e^{-\delta t}
  +2\,\nu \int_{0}^{T}E \big\{|\nabla \mathbf{u}^{n}(t)|^{2}\big\}
  e^{-\delta t}\d t \\
  & \;\leq |\mathbf{u}(0)|^2 + \varepsilon |p(0)|^2
  +\int_{0}^{T} \big[\frac{1}{\delta }|\mathbf{f}(t)|^{2}+
  \Tr(\mathbf{g}^{2}(t))\big]e^{-\delta t}\d t,
\end{split}\right.\end{equation}
for any $0\leq t\leq T.$ Moreover, if we suppose
\begin{equation}
  \mathbf{f}\in \mathrm{L}^{\mathrm{p}}(0,T;\mathbb{L}^2(\mathcal{O})),\;
  \mathbf{g}\in \mathrm{L}^{\mathrm{p}}(0,T;\ell_{2}(\mathbb{L}^2(\mathcal{O})))
\end{equation}
then we also have
\begin{equation}\left\{\begin{split}\label{e2.10}
  E\big\{\sup_{0\leq t\leq T}&\big[|\mathbf{u}^{n}(t)|^{\mathrm{p}}
  + \varepsilon |p^n(t)|^{\mathrm{p}}\big] e^{-\delta t} \\
  &\quad + \mathrm{p}\,\nu \int_{0}^{T}|\nabla \mathbf{u}^{n}(t)|^{2}|
  \mathbf{u}^{n}(t)|^{\mathrm{p}-2}e^{-\delta t}\d t\big\}\leq
  |\mathbf{u}(0)|^{\mathrm{p}}  \\
  &\quad +\varepsilon |p(0)|^{\mathrm{p}}
  + C_{\delta ,\mathrm{p},T}\int_{0}^{T} \big[|\mathbf{f}(t)|^{\mathrm{p}} +
  \Tr(\mathbf{g}^{2}(t))^{\mathrm{p}/2}\big] e^{-\delta t}\d t,
\end{split}\right.\end{equation}
for some constant $C_{\delta ,\mathrm{p},T}$ depending only on
$\delta>0,$ $\varepsilon>0,$ $1\leq \mathrm{p}<\infty $ and $T>0.$
\end{proposition}

\begin{proof}
From \eqref{e2.4} we notice that,
\begin{align}\label{e2.11}
& d(\mathbf{u}^n(t),\mathbf{u}^n(t)) + \langle
-\nu\triangle\mathbf{u}^n(t),\mathbf{u}^n(t)\rangle\d t + \langle
(\mathbf{u}^n(t)\cdot\nabla)\mathbf{u}^n(t),\mathbf{u}^n(t)\rangle\d
t\nonumber\\
& \qquad\qquad + \frac{1}{2}\langle
\mathbf{u}^n(t)\Div\mathbf{u}^n(t),\mathbf{u}^n(t)\rangle\d t +
\langle\nabla p^n(t),\mathbf{u}^n(t)\rangle\d t \nonumber\\
& \qquad \qquad = ( \mathbf{f}(t),\mathbf{u}^n(t) ) \d t + \sum_k (
  \mathbf{g}_k(t),\mathbf{u}^n(t)) \d w_k(t).
\end{align}
It is clear that
$$\langle-\nu\triangle\mathbf{u}^n(t),\mathbf{u}^n(t)\rangle =
\nu|\nabla\mathbf{u}^n(t)|^{2},$$ and the equation \eqref{e1.8}
yields
\begin{equation}
\langle(\mathbf{u}^n(t)\cdot\nabla)\mathbf{u}^n(t),\mathbf{u}^n(t)\rangle
= -\langle\mathbf{u}^n(t)\Div\mathbf{u}^n(t),\mathbf{u}^n(t)\rangle
-\langle(\mathbf{u}^n(t)\cdot\nabla)\mathbf{u}^n(t),\mathbf{u}^n(t)\rangle.\nonumber
\end{equation}
Hence
\begin{align}\label{e2.12}
\langle(\mathbf{u}^n(t)\cdot\nabla)\mathbf{u}^n(t),\mathbf{u}^n(t)\rangle
+ \frac{1}{2}\langle
\mathbf{u}^n(t)\Div\mathbf{u}^n(t),\mathbf{u}^n(t)\rangle = 0.
\end{align}
 Using the equation $\varepsilon\dot{p}^n(t) +
\Div\mathbf{u}^n(t) = 0$ we get from \eqref{e1.6}
\begin{equation}\label{e2.13} \langle\nabla
p^{n}(t),\mathbf{u}^n(t)\rangle = -\langle p^{n}(t),
\Div\mathbf{u}^n(t)\rangle = -\langle
p^{n}(t),-\varepsilon\dot{p}^n(t)\rangle = \frac{\varepsilon}{2}\
\frac{\d}{\d t}|{p}^n(t)|^{2}.
\end{equation}
Combining all the above results one can get from \eqref{e2.11}
\begin{align}
& \frac{1}{2}\ \d |\mathbf{u}^n(t)|^{2} +
\nu|\nabla\mathbf{u}^n(t)|^{2}\d t + \frac{\varepsilon}{2}\
\d|{p}^n(t)|^{2} \nonumber\\
 & \qquad\qquad\qquad = ( \mathbf{f}(t),\mathbf{u}^n(t) ) \d t +
\frac{1}{2} \ \mathrm{Tr}(\mathbf{g}^{2}(t))\d t + \sum_k (
  \mathbf{g}_k(t),\mathbf{u}^n(t)) \d w_k(t).\nonumber
\end{align}
Rearranging the terms we get the desired energy equality
\eqref{e2.7}.

Next, we calculate the stochastic differential of the process $$
\mathrm{F}(t) := \big[|\mathbf{u}^n(t)|^{2} + \varepsilon
|{p}^n(t)|^{2}\big]e^{-\delta t} $$ to get
$$\d \mathrm{F}(t) = e^{-\delta t}\d \big[|\mathbf{u}^n(t)|^{2} + \varepsilon
|{p}^n(t)|^{2}\big] - \delta\big[|\mathbf{u}^n(t)|^{2} + \varepsilon
|{p}^n(t)|^{2}\big]e^{-\delta t}\d t.$$ Using the energy equality
\eqref{e2.7} we have
\begin{align}
\d \mathrm{F}(t)= &-2\nu|\nabla\mathbf{u}^{n}(t)|^{2}e^{-\delta t}\d
t + 2\,(\mathbf{f}(t), \mathbf{u}^{n}(t))e^{-\delta t}\d t +
\Tr(\mathbf{g}^{2}(t))e^{-\delta t}\d t \nonumber\\ & +
2\sum_{k}(\mathbf{g}_{k}(t),\mathbf{u}^{n}(t)) e^{-\delta t} \d
w_{k}(t) - \delta\big[|\mathbf{u}^n(t)|^{2} + \varepsilon
|{p}^n(t)|^{2}\big]e^{-\delta t}\d t.\label{e2.14}
\end{align}
Now using the inequality $$ 2ab \leq \delta a^{2} + \frac{1}{\delta}
b^{2}$$ on $2\,(\mathbf{f}(t), \mathbf{u}^{n}(t))$ we have
\begin{align}
2\,(\mathbf{f}(t), \mathbf{u}^{n}(t)) \leq \delta
|\mathbf{u}^n(t)|^{2} + \frac{1}{\delta} |\mathbf{f}(t)|^{2} \leq
\delta\big[|\mathbf{u}^n(t)|^{2} + \varepsilon |{p}^n(t)|^{2}\big] +
\frac{1}{\delta} |\mathbf{f}(t)|^{2}.\nonumber
\end{align}
Then \eqref{e2.14} yields
\begin{align}
\d \mathrm{F}(t)= &-2\nu|\nabla\mathbf{u}^{n}(t)|^{2}e^{-\delta t}\d
t + \delta\big[|\mathbf{u}^n(t)|^{2} + \varepsilon
|{p}^n(t)|^{2}\big]e^{-\delta t}\d t \nonumber\\ & +
\frac{1}{\delta} |\mathbf{f}(t)|^{2}e^{-\delta t}\d t +
\Tr(\mathbf{g}^{2}(t))e^{-\delta t}\d t +
2\sum_{k}(\mathbf{g}_{k}(t),\mathbf{u}^{n}(t)) e^{-\delta t} \d
w_{k}(t) \nonumber\\ & - \delta\big[|\mathbf{u}^n(t)|^{2} +
\varepsilon |{p}^n(t)|^{2}\big]e^{-\delta t}\d t.\nonumber
\end{align}
Rearranging the terms we have
\begin{align}
& \d \big[\big\{|\mathbf{u}^n(t)|^{2} + \varepsilon
|{p}^n(t)|^{2}\big\}e^{-\delta t}\big] +
2\nu|\nabla\mathbf{u}^{n}(t)|^{2}e^{-\delta t}\d t \nonumber\\
& \qquad \leq\ \frac{1}{\delta} |\mathbf{f}(t)|^{2}e^{-\delta t}\d t
+ \Tr(\mathbf{g}^{2}(t))e^{-\delta t}\d t +
2\sum_{k}(\mathbf{g}_{k}(t),\mathbf{u}^{n}(t)) e^{-\delta t} \d
w_{k}(t). \nonumber
\end{align}
Next we integrate in $[0, T]$ to get
\begin{align}
& \big[|\mathbf{u}^n(t)|^{2} + \varepsilon
|{p}^n(t)|^{2}\big]e^{-\delta t} +
2\nu\int_{0}^{T}|\nabla\mathbf{u}^{n}(t)|^{2}e^{-\delta t}\d t
\nonumber\\
& \qquad \leq\ |\mathbf{u}(0)|^{2} + \varepsilon |{p}(0)|^{2} +
\int_{0}^{T}\big[\frac{1}{\delta} |\mathbf{f}(t)|^{2} +
\Tr(\mathbf{g}^{2}(t))\big]e^{-\delta t}\d t \nonumber\\
& \qquad\quad +
2\sum_{k}\int_{0}^{T}(\mathbf{g}_{k}(t),\mathbf{u}^{n}(t))
e^{-\delta t} \d w_{k}(t). \nonumber
\end{align}
Finally taking mathematical expectation and keeping in mind the
expectation of a stochastic integral is zero, we have the desired
result \eqref{e2.8}. \\
Similarly, consider $$ \mathrm{G}(t) :=
\big[|\mathbf{u}^n(t)|^{\mathrm{p}} + \varepsilon
|{p}^n(t)|^{\mathrm{p}}\big]e^{-\delta t} $$ and use It\^o calculus.
Here we check that its stochastic differential satisfies
\begin{align}
\d \mathrm{G}(t) = &-\delta\big[|\mathbf{u}^n(t)|^{\mathrm{p}} +
\varepsilon
|{p}^n(t)|^{\mathrm{p}}\big]e^{-\delta t}\d t\nonumber\\
&+ \frac{\mathrm{p}}{2}|\mathbf{u}^n(t)|^{\mathrm{p}-2}\big\{\d
\big[|\mathbf{u}^n(t)|^{2} + \varepsilon
|{p}^n(t)|^{2}\big]\big\}e^{-\delta t}\nonumber\\
&+
\frac{\mathrm{p}(\mathrm{p}-1)}{8}|\mathbf{u}^n(t)|^{\mathrm{p}-4}\big\{\d
\big[|\mathbf{u}^n(t)|^{2} + \varepsilon
|{p}^n(t)|^{2}\big]\big\}^{2}e^{-\delta t}.\nonumber
\end{align}
Using the energy equality \eqref{e2.7} we get
\begin{align}
\d \mathrm{G}(t)&+\delta\big[|\mathbf{u}^n(t)|^{\mathrm{p}} +
\varepsilon |{p}^n(t)|^{\mathrm{p}}\big]e^{-\delta t}\d t\nonumber\\
&=\frac{\mathrm{p}}{2}|\mathbf{u}^n(t)|^{\mathrm{p}-2}\Big\{-2\nu \,
  \|\mathbf{u}^{n}(t)\|^{2}\d t
  +\big[2\,(\mathbf{f}(t), \mathbf{u}^{n}(t))+
  \Tr(\mathbf{g}^{2}(t))\big]\d t\nonumber\\
&\quad+ 2\sum_{k}(\mathbf{g}_{k}(t),\mathbf{u}^{n}(t)) \d
w_{k}(t)\Big\}e^{-\delta t}\nonumber\\
&\quad+\frac{\mathrm{p}(\mathrm{p}-1)}{8}|\mathbf{u}^n(t)|^{\mathrm{p}-4}
\big[4\sum_{k}(\mathbf{g}_{k}(t),\mathbf{u}^{n}(t))^{2}\d
t\big]e^{-\delta t}.\nonumber
\end{align}
Simplification and rearrangement of the terms in the above equation
yields
\begin{align}\label{e2.15}
\d \mathrm{G}(t)&+
\nu\mathrm{p}\|\mathbf{u}^{n}(t)\|^{2}|\mathbf{u}^n(t)|^{\mathrm{p}-2}e^{-\delta
t}\d t +\delta\big[|\mathbf{u}^n(t)|^{\mathrm{p}} +
\varepsilon |{p}^n(t)|^{\mathrm{p}}\big]e^{-\delta t}\d t\nonumber\\
&=|\mathbf{u}^n(t)|^{\mathrm{p}-2}\big[\mathrm{p}\,(\mathbf{f}(t),
\mathbf{u}^{n}(t))+\frac{\mathrm{p}^{2}}{2}\Tr(\mathbf{g}^{2}(t))\big]e^{-\delta
t}\d t\nonumber\\
&\quad+\mathrm{p}\sum_{k}(\mathbf{g}_{k}(t),\mathbf{u}^{n}(t))
|\mathbf{u}^n(t)|^{\mathrm{p}-2}e^{-\delta t}\d w_{k}(t).
\end{align}
Now in the first and the second terms on the right hand side of the
above equation we apply the following elementary inequality
$$\lambda \mathrm{a}\mathrm{b} \leq \frac{(\alpha\mathrm{a})^{\mathrm{p}}}{\mathrm{p}} + \frac{(\beta\mathrm{b})
^{\mathrm{q}}}{\mathrm{q}},$$ where $ \qquad \frac{1}{\mathrm{p}} +
\frac{1}{\mathrm{q}} = 1, \qquad \lambda = \alpha\beta > 0, \qquad
\mathrm{a}\mathrm{b} > 0.$ \\
Then choosing $$\lambda = \mathrm{p},\quad \alpha =
\frac{\mathrm{p}}{(\frac{\delta\mathrm{q}}{2})^{^{\frac{1}{\mathrm{q}}}}},
\quad \beta =
\Big(\frac{\delta\mathrm{q}}{2}\Big)^{^{\frac{1}{\mathrm{q}}}}.$$ we
get
\begin{align}
\mathrm{p}\,(\mathbf{f}(t),
\mathbf{u}^{n}(t))|\mathbf{u}^n(t)|^{\mathrm{p}-2} &\leq
\alpha^{\mathrm{p}}\frac{|\mathbf{f}(t)|^{\mathrm{p}}}{\mathrm{p}} +
\beta^{\mathrm{q}}\frac{|\mathbf{u}^n(t)|^{(\mathrm{p}-1)\mathrm{q}}}{\mathrm{q}}\nonumber\\
&= C_{\delta, \mathrm{p}}\ |\mathbf{f}(t)|^{\mathrm{p}} +
\frac{\delta}{2}|\mathbf{u}^n(t)|^{\mathrm{p}}\nonumber,
\end{align}
where the constant $C_{\delta, \mathrm{p}} > 0$ depends only on
$\delta > 0$ and $1\leq\mathrm{p}<\infty$.\\
Similarly with proper choices of $\alpha$ and $\beta$, we can prove
that
\begin{align}
\frac{\mathrm{p}^{2}}{2}\Tr(\mathbf{g}^{2}(t))|\mathbf{u}^n(t)|^{\mathrm{p}-2}
\leq C_{\delta, \mathrm{p}}\
\Tr(\mathbf{g}^{2}(t))^{\frac{\mathrm{p}}{2}} +
\frac{\delta}{2}|\mathbf{u}^n(t)|^{\mathrm{p}}.\nonumber
\end{align}
Then \eqref{e2.15} yields
\begin{align}\label{e2.16}
\d \mathrm{G}(t)&+
\nu\mathrm{p}\|\mathbf{u}^{n}(t)\|^{2}|\mathbf{u}^n(t)|^{\mathrm{p}-2}e^{-\delta
t}\d t\nonumber\\
&\leq C_{\delta, \mathrm{p}}\big[|\mathbf{f}(t)|^{\mathrm{p}} +
\Tr(\mathbf{g}^{2}(t))^{\frac{\mathrm{p}}{2}}\big]e^{-\delta t}\d
t\nonumber\\
&\quad+\mathrm{p}\sum_{k}(\mathbf{g}_{k}(t),\mathbf{u}^{n}(t))
|\mathbf{u}^n(t)|^{\mathrm{p}-2}e^{-\delta t}\d w_{k}(t).
\end{align}
Integrating the stochastic differential \eqref{e2.16},then taking
the sup norm in $[0, T]$ and finally taking the mathematical
expectation we have
\begin{align}\label{e2.17}
 E\Big\{\sup_{0\leq t\leq T}&\big[|\mathbf{u}^{n}(t)|^{\mathrm{p}}
  + \varepsilon |p^n(t)|^{\mathrm{p}}\big] e^{-\delta t}\nonumber \\
  &\quad + \mathrm{p}\,\nu \int_{0}^{T}|\nabla \mathbf{u}^{n}(t)|^{2}|
  \mathbf{u}^{n}(t)|^{\mathrm{p}-2}e^{-\delta t}\d t\Big\}\leq
  |\mathbf{u}(0)|^{\mathrm{p}}\nonumber  \\
  &\quad +\varepsilon |p(0)|^{\mathrm{p}}
  + C_{\delta ,\mathrm{p},T}\int_{0}^{T} \big[|\mathbf{f}(t)|^{\mathrm{p}} +
  \Tr(\mathbf{g}^{2}(t))^{\mathrm{p}/2}\big] e^{-\delta t}\d
  t\nonumber\\
  &\quad + \mathrm{p}\ E\Big\{\sup_{0\leq t\leq T}\Big|\int_{0}^{t}\sum_{k}(\mathbf{g}_{k}(s),\mathbf{u}^{n}(s))
|\mathbf{u}^n(s)|^{\mathrm{p}-2}e^{-\delta s}\d w_{k}(s)\Big|\Big\}.
\end{align}
By means of martingale inequality, we deduce
\begin{align}
E\Big\{\sup_{0\leq t\leq
T}&\Big|\int_{0}^{t}\sum_{k}(\mathbf{g}_{k}(s),\mathbf{u}^{n}(s))
|\mathbf{u}^n(s)|^{\mathrm{p}-2}e^{-\delta s}\d
w_{k}(s)\Big|\Big\}\nonumber\\
&\quad\leq C\
E\Big\{\Big(\int_{0}^{T}\sum_{k}\big[(\mathbf{g}_{k}(t),\mathbf{u}^{n}(t))
|\mathbf{u}^n(t)|^{\mathrm{p}-2}e^{-\delta t}\big]^{2}\d
t\Big)^{1/2}\Big\}\nonumber\\
&\quad\leq C\
E\Big\{\Big(\int_{0}^{T}\Tr(\mathbf{g}^{2}(t))|\mathbf{u}^n(t)|^{\mathrm{2p}-2}e^{-2\delta
t}\d t\Big)^{1/2}\Big\}\nonumber\\
&\quad\leq C\ E\Big\{\sup_{0\leq t\leq T}
(|\mathbf{u}^n(t)|^{\mathrm{p}-1}e^{-\delta
t/\mathrm{p^{\prime}}})\Big(\int_{0}^{T}\Tr(\mathbf{g}^{2}(t))e^{-2\delta
t/\mathrm{p}}\d t\Big)^{1/2}\Big\}\nonumber\\
&\quad\leq \frac{\delta}{2}\ E\big\{\sup_{0\leq t\leq T}
(|\mathbf{u}^n(t)|^{\mathrm{p}}e^{-\delta t})\big\} + C_{\delta
,\mathrm{p},T}\
E\big\{\int_{0}^{T}\Tr(\mathbf{g}^{2}(t))^{\mathrm{p}/2}e^{-\delta
t}\d t\big\},\label{e2.18}
\end{align}
where the constant $C_{\delta ,\mathrm{p},T}$ depends only on
$\delta>0,$ $1\leq \mathrm{p}<\infty $ and $T>0.$ \\
Using \eqref{e2.18} in \eqref{e2.17} we get the desired estimate
\eqref{e2.10}.
\end{proof}

Now we deal with the existence and uniqueness of the SPDE and its
finite-dimensional approximation.

\begin{proposition}[uniqueness]
Let $\mathbf{u}$ be a solution of the stochastic Navier-Stokes
equation (SPDE) \eqref{e2.2} with the regularity
\begin{equation}\left\{\begin{split}
  &\mathbf{u}\in \mathrm{L}^2(\Omega ;C^{0}(0,T;\mathbb{L}^2(\mathcal{O}))
  \cap\mathrm{L}^2(0,T;\mathbb{H}^1_0(\mathcal{O}))),\quad \\
  &\mathbf{u}\in \mathbb{L}^4(\Omega\times\mathcal{O}\times(0,T)),
  \quad p\in \mathrm{L}^2(\Omega\times\mathcal{O}\times (0,T)),\;
\end{split}\right.\end{equation}
and let the data $\mathbf{f},$ $\mathbf{g},$ $\mathbf{u}_{0}$ and
$p_0$ satisfy the condition
\begin{equation}\left\{\begin{split}
  &\mathbf{f}\in
  \mathrm{L}^2(0,T;\mathbb{H}^{-1}(\mathcal{O})),\quad
  \mathbf{g} \in \mathrm{L}^2(0,T;\ell
  _{2}(\mathbb{L}^2(\mathcal{O}))),\\
  &\mathbf{u}_{0}\in \mathbb{L}^2(\mathcal{O}), \quad
  p_0\in \mathrm{L}^2(\mathcal{O}).
\end{split}\right.\end{equation}
If $\mathbf{v}$ in $\mathrm{L}^2(\Omega;C^0(0,T,\mathbb{L}^2(\mathcal{O}))\cap %
\mathrm{L}^2(0,T,\mathbb{H}^1_0(\mathcal{O})))$ is another solution
of the stochastic Navier-Stokes equation \eqref{e2.2}, then
\begin{equation}\left\{\begin{split}\label{e2.21}
  \big[|\mathbf{u}(t)-\mathbf{v}(t)|^2+\varepsilon|p(t)-q(t)|^2 \big]\,\exp
  \Big[-\frac{27}{\nu ^3}
  \int_{0}^{t}\Vert \mathbf{u}(s)\Vert _{_{\mathbb{L}^{4}(\mathcal{O})}}^{4}
  \d s \Big]\leq \\
  \leq |\mathbf{u}(0)-\mathbf{v}(0)|^2 + \varepsilon|p(0)-q(0)|^2,
\end{split}\right.\end{equation}
with probability $1$ \ for any $0\leq t\leq T$ and $\varepsilon>0.$
\end{proposition}

\begin{proof}
 Indeed  if $\mathbf{u}$ and
$\mathbf{v}$ are two solutions then $\mathbf{w}=
\mathbf{v}-\mathbf{u}$ solves the deterministic equation
\[
  \partial_t \mathbf{w}(t) - \nu\triangle \mathbf{w}(t)
  + \nabla(q(t) - p(t)) = \hat{B}(\mathbf{u})
  - \hat{B}(\mathbf{v}) \quad \text{in }
  \mathbb{L}^2(0,T;\mathbb{H}^{-1}(\mathcal{O})),
\]
with $\hat{B}(\mathbf{u})=[(\mathbf{u}\cdot\nabla) + \frac{1}{2}
\Div\mathbf{u}]\mathbf{u}.$  Notice that actually $p$ and $q$ are
better processes, they belong to the space
$\mathrm{L}^2(\Omega;C^{0}(0,T;\mathbb{L}^2(\mathcal{O}))).$

Next, setting $$r(t):=\frac{27}{\nu^3}\int_0^t \|\mathbf{u(s)}
\|_{_{\mathbb{L}^4(\mathcal{O})}}^4 \d s$$ we have
\begin{align}
\d \langle\mathbf{w}(t), 2e^{-r(t)}&\mathbf{w}(t)\rangle -
\nu\langle\triangle \mathbf{w}(t), 2e^{-r(t)}\mathbf{w}(t)\rangle\d
t \nonumber\\ &\qquad\qquad+ \langle\nabla(q(t) - p(t)),
2e^{-r(t)}\mathbf{w}(t)\rangle\d
t\nonumber\\
&\qquad\qquad\quad = \langle\hat{B}(\mathbf{u})
  - \hat{B}(\mathbf{v}), 2e^{-r(t)}\mathbf{w}(t)\rangle\d
t.\nonumber
\end{align}
Using \eqref{e1.6} we get
\begin{align}\label{e2.22}
e^{-r(t)}\d (|\mathbf{w}(t)|^{2}) + & 2\nu
e^{-r(t)}\|\mathbf{w}(t)\|^{2}\d t + 2e^{-r(t)}\langle p(t) - q(t),
\Div\mathbf{w}(t)\rangle\d t\nonumber\\
&\qquad\qquad = 2e^{-r(t)}\langle\hat{B}(\mathbf{u})
  - \hat{B}(\mathbf{v}), \mathbf{w}(t)\rangle\d
t.
\end{align}
Since $\varepsilon\partial_t(q(t) - p(t)) + \Div\mathbf{w}(t) = 0$,
Lemma 2.3 and \eqref{e2.22} yield
\begin{align}
& e^{-r(t)}\d \Big[|\mathbf{w}(t)|^{2} + \varepsilon |q(t) -
p(t)|^{2}\Big]\nonumber\\
&\qquad = -2\nu e^{-r(t)}\|\mathbf{w}(t)\|^{2}\d t -
2e^{-r(t)}\langle\hat{B}(\mathbf{v})
  - \hat{B}(\mathbf{u}), \mathbf{w}(t)\rangle\d
t\nonumber\\
&\qquad \leq -2\nu e^{-r(t)}\|\mathbf{w}(t)\|^{2}\d t +
2e^{-r(t)}\Big[\frac{\nu}{2}\|\mathbf{w}(t)\|^{2} +
\frac{27}{2\nu^{3}}|\mathbf{w}(t)|^{2}\|\mathbf{u(t)}
\|_{_{\mathbb{L}^4(\mathcal{O})}}^4\Big]\d t\nonumber\\
&\qquad = -\nu e^{-r(t)}\|\mathbf{w}(t)\|^{2}\d t +
\dot{r}(t)e^{-r(t)}|\mathbf{w}(t)|^{2}\d t\nonumber\\
&\qquad\leq-\nu e^{-r(t)}\|\mathbf{w}(t)\|^{2}\d t +
\dot{r}(t)e^{-r(t)}\Big[|\mathbf{w}(t)|^{2} + \varepsilon |q(t) -
p(t)|^{2}\Big]\d t.\nonumber
\end{align}
Hence
\begin{align}
\d \Big[e^{-r(t)}\Big\{|\mathbf{w}(t)|^{2} + \varepsilon |q(t) -
p(t)|^{2}\Big\}\Big] \leq 0.\nonumber
\end{align}
Hence, integrating in $t$, we deduce \eqref{e2.21}, with probability
$1$.
\end{proof}

Each solution $\mathbf{u}$ in the space $\mathrm{L}^2(\Omega
;\mathrm{L}^{\infty }(0,T;\mathbb{H}^{-1}(\mathcal{O}))\cap
\mathrm{L}^2(0,T;\mathbb{H}^1_0(\mathcal{O})))$ of the stochastic
Navier-Stokes equation actually belongs to a better space, namely
the space
$\mathrm{L}^2(\Omega;C^{0}(0,T;\mathbb{L}^2(\mathcal{O}))\cap
\mathbb{L}^{4}(\mathcal{O}\times (0,T)))$ in 2-D,
$\mathcal{O}\subset \mathbb{R}^{2}.$ Thus in 2-D, the uniqueness
holds in the space
$\mathrm{L}^2(\Omega;\mathrm{L}^2(0,T;\mathbb{H}^{-1}(\mathcal{O}))).$

If a given adapted process $\mathbf{u}$ in
$\mathrm{L}^2(\Omega;\mathrm{L}^{\infty
}(0,T;\mathbb{L}^{2}(\mathcal{O}))\cap
\mathrm{L}^2(0,T;\mathbb{H}^1_0(\mathcal{O})))$ satisfies
\begin{equation}
  \d(\mathbf{u}(t),\mathbf{v})=\langle \mathbf{F}(t),\mathbf{v}\rangle
  \d t +(\mathbf{g}(t),\mathbf{v})\d w(t),
\end{equation}
for any function $\mathbf{v}$ in $\mathbb{H}^1_0(\mathcal{O})$ and
some functions $\mathbf{F}$ in
$\mathrm{L}^2(0,T;\mathbb{H}^{-1}(\mathcal{O}))$ and $\mathbf{g}$ in
$\mathrm{L}^2(0,T;\ell _{2}(\mathbb{L}^2(\mathcal{O}))),$ then we
can find a version of $\mathbf{u}$ (which is still denoted by
$\mathbf{u}$) in $\mathrm{L}^2(\Omega
;C^0(0,T;\mathbb{L}^2(\mathcal{O})))$ satisfying the energy equality
\begin{equation}
  \d|\mathbf{u}(t)|^{2}=\big[2\langle \mathbf{F}(t),\mathbf{u}(t)\rangle +
  \Tr(\mathbf{g^{2}}(t)\big]\d t+
  2(\mathbf{g}(t),\mathbf{u}(t))\d w(t)
\end{equation}
see e.g. Gyongy and Krylov~\cite{GyKr82}, Pardoux~\cite{Pa79}.

\begin{definition}($Strong\ Solution$)
A strong solution $\mathbf{u}$ is defined on a given probability
space ($\Omega$, $\Sigma$, $\Sigma_{t}$, $\mathcal{M}$) as a
$\mathrm{L}^2(\Omega;\mathrm{L}^{\infty
}(0,T;\mathbb{L}^{2}(\mathcal{O}))\cap
\mathrm{L}^2(0,T;\mathbb{H}^1_0(\mathcal{O}))\cap
C^0(0,T;\mathbb{L}^{2}(\mathcal{O})))$ valued function which
satisfies the stochastic Navier-Stokes equation \eqref{e2.2} in the
weak sense and also the energy inequality
\begin{align*}
  &E\big\{\sup_{0\leq t\leq T}\big[|\mathbf{u}(t)|^{\mathrm{p}}
  + \varepsilon |p(t)|^{\mathrm{p}}\big] e^{-\delta t} + \mathrm{p}\,\nu \int_{0}^{T}|\nabla \mathbf{u}(t)|^{2}|
  \mathbf{u}(t)|^{\mathrm{p}-2}e^{-\delta t}\d t\big\}\\
  &\quad\leq
  |\mathbf{u}(0)|^{\mathrm{p}} +\varepsilon |p(0)|^{\mathrm{p}}
  + C_{\delta ,\mathrm{p},T}\int_{0}^{T} \big[|\mathbf{f}(t)|^{\mathrm{p}} +
  \Tr(\mathbf{g}^{2}(t))^{\mathrm{p}/2}\big] e^{-\delta t}\d t,
\end{align*}
where the constant $C_{\delta ,\mathrm{p},T}$ depends only on
$\delta>0,$ $\varepsilon>0,$ $1\leq \mathrm{p}<\infty $ and $T>0.$
\end{definition}


\begin{proposition}[2-D existence]
Let $\mathbf{f},$ $\mathbf{g}$ and $\mathbf{u}_{0}$ be such that
\begin{equation}\left\{\begin{split}
  &\mathbf{f}\in \mathrm{L}^{\mathrm{p}}(0,T;\mathbb{H}^{-1}(\mathcal{O})),
  \quad \mathbf{g} \in \mathrm{L}^{\mathrm{p}}(0,T;\ell _{2}
  (\mathbb{L}^2(\mathcal{O}))), \\
  &\mathbf{u}_{0}\in \mathbb{L}^2(\mathcal{O}), \quad
  p_0\in \mathrm{L}^2(\mathcal{O}).
\end{split}\right.\end{equation}
for some $\mathrm{p}\geq 4.$ Then there is adapted processes
$\mathbf{u}(t,x,\omega)$ and $p(t,x,\omega)$ with the regularity
\begin{equation}\left\{\begin{split}
  &\mathbf{u}\in \mathrm{L}^{\mathrm{p}}(\Omega ;
  C^{0}(0,T;\mathbb{L}^2(\mathcal{O}))) \cap\mathrm{L}^{2}
  (\Omega;\mathrm{L}^{2}(0,T;\mathbb{H}^1_0(\mathcal{O}))),\\
  &p,\;\dot{p}\in \mathrm{L}^{2}(\Omega;\mathrm{L}^2(0,T;
  \mathrm{L}^2(\mathcal{O})))
\end{split}\right.\end{equation}
satisfying the stochastic Navier-Stokes equation \eqref{e2.2} and
the a priori bound \eqref{e2.10} for every $\varepsilon>0.$
\end{proposition}

\begin{proof}
Denoting
$$F(\mathbf{u}):=A\mathbf{u}+\hat{B}(\mathbf{u})-\mathbf{f}:=
-\nu\triangle\mathbf{u}+[(\mathbf{u}\cdot\nabla)+(1/2)\Div\mathbf{u}]\mathbf{u}-\mathbf{f}$$
we have
$$\d \mathbf{u}^{n}(t) + F(\mathbf{u}^{n}(t))\d t + \nabla
p^{n}(t)\d t = \mathbf{g}(t)\d w(t).$$ Then using the a priori
estimate \eqref{e2.10}, it follows from the Banach-Alaoglu theorem
that along a subsequence, the Galerkin approximations
$\{\mathbf{u}^{n}\}$ have the following limits:\\
\begin{align}
&\mathbf{u}^{n}\longrightarrow \mathbf{u}\quad  \text {weakly star
in}\ \mathrm{L}^{\mathrm{p}}(\Omega ;
  \mathrm{L}^{\infty}(0,T;\mathbb{L}^2(\mathcal{O}))) \cap\mathrm{L}^{2}
  (\Omega;\mathrm{L}^{2}(0,T;\mathbb{H}^1_0(\mathcal{O}))),\nonumber\\
& p^{n}\longrightarrow p\ \text{weakly
in}\ \mathrm{L}^{2}(\Omega;\mathrm{L}^{2}(0,T;\mathbb{L}^{2}(\mathcal{O}))),\nonumber\\
& F(\mathbf{u}^{n})\longrightarrow F_{0}\quad \text{weakly in}\
\mathrm{L}^{2}(\Omega;\mathrm{L}^{2}(0,T;\mathbb{H}^{-1}(\mathcal{O}))),
\nonumber
\end{align}
where $\mathbf{u}$ has the It\^o differential
\begin{align}
\d \mathbf{u}(t) + F_{0}(t)\d t + \nabla p(t)\d t = \mathbf{g}(t)\d
w(t)\qquad \text{in}\
\mathrm{L}^{2}(\Omega;\mathrm{L}^{2}(0,T;\mathbb{H}^{-1}(\mathcal{O})))\nonumber
\end{align}
and the energy equality holds, i.e.,
\begin{align}
\d \big[|\mathbf{u}(t)|^{2}+\varepsilon |p(t)|^{2}\big] + 2\langle
F_{0}(t), \mathbf{u}(t)\rangle\d t = \Tr(\mathbf{g}^{2}(t))\d t + 2
(\mathbf{g}(t), \mathbf{u}(t))\d w(t).\nonumber
\end{align}
Now, for any adapted process $\mathbf{v}(t, x, \omega)$ in
$\mathrm{L}^{\infty}((0, T)\times\Omega;
\mathbb{L}^{2}(\mathcal{O})),$ we define
$$r(t,\omega):=\frac{27}{\nu^3}\int_0^t \|\mathbf{v(s,.,\omega)}
\|_{_{\mathbb{L}^4(\mathcal{O})}}^4 \d s$$ as an adapted, continuous
(and bounded in $\omega$) real-valued process in $[0, T]$. Then from
the energy equality
\begin{align}
\d \Big[e^{-r(t)}\big\{|\mathbf{u}^{n}(t)|^{2}+&\varepsilon
|p^{n}(t)|^{2}\big\}\Big] + e^{-r(t)}\langle
2F(\mathbf{u}^{n}(t))+\dot{r}(t)\mathbf{u}^{n}(t),
\mathbf{u}^{n}(t)\rangle\d t\nonumber\\
&+ \varepsilon |p^{n}(t)|^{2}\dot{r}(t)e^{-r(t)}\d t\nonumber\\
&\quad=\Tr(\mathbf{g}^{2}(t))e^{-r(t)}\d t + 2 (\mathbf{g}(t),
\mathbf{u}^{n}(t))e^{-r(t)}\d w(t).\nonumber
\end{align}
Integrating between $0\leq t\leq T$ and taking the mathematical
expectation we have
\begin{align}
&E\Big[e^{-r(T)}\big\{|\mathbf{u}^{n}(T)|^{2}+\varepsilon
|p^{n}(T)|^{2}\big\}-|\mathbf{u}^{n}(0)|^{2}-\varepsilon
|p^{n}(0)|^{2}\Big]\nonumber\\
&\qquad\qquad + E\Big[\int_{0}^{T}e^{-r(t)}\langle
2F(\mathbf{u}^{n}(t))+\dot{r}(t)\mathbf{u}^{n}(t),
\mathbf{u}^{n}(t)\rangle\d t\Big] \nonumber\\
&\qquad\qquad+
E\Big[\varepsilon\int_{0}^{T}|p^{n}(t)|^{2}\dot{r}(t)e^{-r(t)}\d
t\Big]\nonumber\\
&\qquad\qquad\quad =
E\Big[\int_{0}^{T}\Tr(\mathbf{g}^{2}(t))e^{-r(t)}\d t\Big].\nonumber
\end{align}
Considering the fact that the initial conditions $\mathbf{u}^{n}(0)$
and $p^{n}(0)$ converge to $\mathbf{u}(0)$ and $p(0)$ respectively
in $\mathbb{L}^{2}$, and the lower-semi-continuity of the
$\mathbb{L}^{2}$-norm, we deduce
\begin{align}
&\lim_{n}\inf E\Big[-\int_{0}^{T}e^{-r(t)}\langle
2F(\mathbf{u}^{n}(t))+\dot{r}(t)\mathbf{u}^{n}(t),
\mathbf{u}^{n}(t)\rangle\d t\Big] \nonumber\\
&\quad=\lim_{n}\inf
E\Big[e^{-r(T)}\big\{|\mathbf{u}^{n}(T)|^{2}+\varepsilon
|p^{n}(T)|^{2}\big\}-|\mathbf{u}^{n}(0)|^{2}-\varepsilon
|p^{n}(0)|^{2}\nonumber\\
&\qquad\qquad\qquad\qquad +
\varepsilon\int_{0}^{T}|p^{n}(t)|^{2}\dot{r}(t)e^{-r(t)}\d t -
\int_{0}^{T}\Tr(\mathbf{g}^{2}(t))e^{-r(t)}\d t\Big]\nonumber\\
&\quad\geq E\Big[e^{-r(T)}\big\{|\mathbf{u}(T)|^{2}+\varepsilon
|p(T)|^{2}\big\}-|\mathbf{u}(0)|^{2}-\varepsilon
|p(0)|^{2}\nonumber\\
&\qquad\qquad\qquad\qquad +
\varepsilon\int_{0}^{T}|p(t)|^{2}\dot{r}(t)e^{-r(t)}\d t -
\int_{0}^{T}\Tr(\mathbf{g}^{2}(t))e^{-r(t)}\d t\Big]\nonumber\\
&\quad=E\Big[-\int_{0}^{T}e^{-r(t)}\langle
2F_{0}(t)+\dot{r}(t)\mathbf{u}(t), \mathbf{u}(t)\rangle\d
t\Big]\nonumber
\end{align}
Next, by monotonicity on $\mathbb{L}^{4}$-balls, i.e. by Lemma 2.4,
we have
\begin{align}
&\ 2E\Big[\int_{0}^{T}e^{-r(t)}\langle
F(\mathbf{u}^{n}(t))-F(\mathbf{v}(t)),
\mathbf{u}^{n}(t)-\mathbf{v}(t)\rangle \d t\Big]\nonumber\\
&\qquad\qquad\qquad\qquad +
E\Big[e^{-r(t)}\dot{r}(t)|\mathbf{u}^{n}(t)-\mathbf{v}(t)|^{2}\d
t\Big]\ \geq 0.\nonumber
\end{align}
Rearranging the terms we find
\begin{align}
& E\Big[\int_{0}^{T}e^{-r(t)}\langle
2F(\mathbf{v}(t))+\dot{r}(t)\mathbf{v}(t),
\mathbf{v}(t)-\mathbf{u}^{n}(t)\rangle\d t\Big]\nonumber\\
&\quad\geq E\Big[\int_{0}^{T}e^{-r(t)}\langle
2F(\mathbf{u}^{n}(t))+\dot{r}(t)\mathbf{u}^{n}(t),
\mathbf{v}(t)-\mathbf{u}^{n}(t)\rangle\d t\Big]\nonumber
\end{align}
Taking limit in $n$ , we get
\begin{align}
& E\Big[\int_{0}^{T}e^{-r(t)}\langle
2F(\mathbf{v}(t))+\dot{r}(t)\mathbf{v}(t),
\mathbf{v}(t)-\mathbf{u}(t)\rangle\d t\Big]\nonumber\\
&\quad\geq E\Big[\int_{0}^{T}e^{-r(t)}\langle
2F_{0}(t)+\dot{r}(t)\mathbf{u}(t),
\mathbf{v}(t)-\mathbf{u}(t)\rangle\d t\Big]\nonumber
\end{align}
Now we take $\mathbf{v}:=\mathbf{u}+\lambda\mathbf{w}$ with $\lambda
> 0$ and $\mathbf{w}$ is an adapted process in\\
$\mathrm{L}^{4}(\Omega;
\mathrm{L}^{\infty}(0,T;\mathbb{L}^{2}(\mathcal{O})))\cap\mathrm{L}^{2}(\Omega;
\mathrm{L}^{2}(0,T; \mathbb{H}_{0}^{1}(\mathcal{O})))$.\\
Then we have
\begin{align}
&\lambda E\Big[\int_{0}^{T}e^{-r(t)}\langle
2F\big(\mathbf{u}(t)+\lambda\mathbf{w}(t)\big)-2F_{0}(t),
\mathbf{w}(t)\rangle\d t\Big]\nonumber\\
&\qquad\qquad\qquad\qquad +
\lambda^{2}E\Big[\int_{0}^{T}e^{-r(t)}\dot{r}(t)|\mathbf{w}(t)|^{2}\d
t\Big]\ \geq 0.\nonumber
\end{align}
Dividing by $\lambda$ on both sides of the inequality above, and
letting $\lambda$ go to $0$, one obtains
\begin{align}
E\Big[\int_{0}^{T}e^{-r(t)}\langle F(\mathbf{u}(t))-F_{0}(t),
\mathbf{w}(t)\rangle\d t\Big]\ \geq 0.\nonumber
\end{align}
Since $\mathbf{w}$ is arbitrary, we conclude that
$F_{0}(t)=F(\mathbf{u}(t))$. Thus the existence of a strong solution
of the stochastic Navier-Stokes equation \eqref{e2.2} has been
proved.
\end{proof}

\section{Convergence as $\varepsilon\rightarrow 0$}

 We will now study the asymptotic limit of $\varepsilon\rightarrow
 0$. Let us consider the family of perturbed systems(depending on
 the positive parameter $\varepsilon$).
\begin{align}\label{e2.28}
  &\partial _{t}\mathbf{u^{\varepsilon}}-\nu\triangle\mathbf{u^{\varepsilon}}+
  \sum_{i=1}^{2}\ u_{i}^{\varepsilon}D_{i}\mathbf{u^{\varepsilon}}
  +\frac{1}{2}(\Div\mathbf{u^{\varepsilon}})\mathbf{u^{\varepsilon}}
  +\nabla p^{\varepsilon} =\mathbf{f}+\mathbf{g}(t)\d w(t)\;\\
  &\qquad\qquad\qquad\qquad\qquad\qquad\qquad\qquad\qquad\qquad \text{ in }\;
  \mathrm{L}^2(0,T;\mathbb{H}^{-1}(\mathcal{O})),\nonumber\\
  &\varepsilon\partial_{t}p^{\varepsilon} + \Div\mathbf{u^{\varepsilon}}=0 \;
  \text{ in }\; \mathrm{L}^2(0,T;\mathrm{L}^2(\mathcal{O})),\label{e2.29}
\end{align} with the initial
conditions
\begin{equation}\label{e2.30}
  \mathbf{u^{\varepsilon}}(0)=\mathbf{u}_{0}\quad
  \text{in}\quad\mathbb{L}^2(\mathcal{O})\qquad\text{and}\qquad
  p^{\varepsilon}(0)=p_0\quad\text{in}\quad
  \mathrm{L}^2(\mathcal{O}).
\end{equation}
This is a method to overcome the computational difficulties
connected with the constraint "$\Div\mathbf{u} = 0$". The equations
\eqref{e2.28}-\eqref{e2.29} are easier to approximate than the
original stochastic Navier-Stokes equation as the constraint
"$\Div\mathbf{u} = 0$" has been replaced
by the evolution equation \eqref{e2.29}.\\
Here we will show how the solutions of the perturbed problems
converge to the solutions of the incompressible stochastic
Navier-Stokes equation as $\varepsilon\rightarrow 0$. The idea of
the proof is similar to the deterministic case presented in
Temam~\cite{Te95}.\\
Let
\begin{equation}
  \mathbf{f}\in \mathrm{L}^2(0,T;\mathbb{L}^2(\mathcal{O})),\;
  \mathbf{g}\in \mathrm{L}^2(0,T;\ell_{2}(\mathbb{L}^2(\mathcal{O}))),
  \; \mathbf{u}_{0}\in \mathbb{L}^2(\mathcal{O}),\;
  p_0\in \mathrm{L}^2(\mathcal{O}).\!\!\!\nonumber
\end{equation}
Then we can write the above mentioned perturbed systems as It\^o
Stochastic equations in variational form
\begin{equation}
  \left\{\begin{split}\label{e2.31}
  &\d(\mathbf{u^{\varepsilon}}(t),\mathbf{v}) + \langle -\nu\triangle\mathbf{u^{\varepsilon}}(t),
  \mathbf{v} \rangle \, \d t +
  \hat{B}(\mathbf{u^{\varepsilon}}, \mathbf{u^{\varepsilon}}, \mathbf{v})\d t\\
  &\qquad\qquad\qquad\qquad +\langle \nabla p^{\varepsilon}(t), \mathbf{v} \rangle \, \d t
  = ( \mathbf{f},\mathbf{v} ) \, \d t + (\mathbf{g}(t),\mathbf{v})\, \d w(t),\!\!\\
  &\langle \varepsilon \dot{p^{\varepsilon}}(t)+ \Div\mathbf{u^{\varepsilon}}(t),q\rangle = 0,
  \end{split}\right.
\end{equation}
in $(0,T),$ with the initial conditions
\begin{equation}\label{e2.32}
  (\mathbf{u^{\varepsilon}}(0),\mathbf{v}) =
  (\mathbf{u}_0,\mathbf{v})\quad\text{and}\quad
  (p^{\varepsilon}(0),q)=(p_0,q),
\end{equation}
for any $\mathbf{v}$ in the space $\mathbb{H}^1_0(\mathcal{O})$ and
any $q$ in
$\mathrm{L}^{2}(\Omega;\mathrm{L}^{2}(0,T;\mathbb{L}^{2}(\mathcal{O}))).$


\begin{proposition}
As  $\varepsilon\rightarrow 0$, the solutions
$\{\mathbf{u^{\varepsilon}},p^{\varepsilon}\}$ of the equations
\eqref{e2.31}-\eqref{e2.32} converge to the solution $\mathbf{u}$ of
the incompressible stochastic Navier-Stokes equation.
\end{proposition}

\begin{proof}
First we should point out that the solutions
$\{\mathbf{u^{\varepsilon}},p^{\varepsilon}\}$ of the equations
\eqref{e2.31}-\eqref{e2.32} satisfy the monotonicity property in
Lemma 2.4, the energy equality \eqref{e2.7} and the a priori
estimates \eqref{e2.8} and \eqref{e2.10}. By virtue of these a
priori estimates and using the Banach-Alaoglu theorem, along a
subsequence the approximations
$\{\mathbf{u^{\varepsilon}},p^{\varepsilon}\}$ have the following
limits:
\begin{align}\label{e2.33}
&\mathbf{u^{\varepsilon}}\longrightarrow \mathbf{u}\quad  \text
{weakly star in}\ \mathrm{L}^{\mathrm{p}}(\Omega ;
  \mathrm{L}^{\infty}(0,T;\mathbb{L}^2(\mathcal{O}))),\\
&\qquad\qquad\qquad\qquad\qquad\text{weakly in}\ \mathrm{L}^{2}
  (\Omega;\mathrm{L}^{2}(0,T;\mathbb{H}^1_0(\mathcal{O}))),\nonumber\\
& \sqrt{\varepsilon} p^{\varepsilon}\longrightarrow \chi\
\text{weakly in}\
\mathrm{L}^{2}(\Omega;\mathrm{L}^{2}(0,T;\mathbb{L}^{2}(\mathcal{O}))).\label{e2.34}
\end{align}
Let us denote
$$F(\mathbf{u^{\varepsilon}}):=A\mathbf{u^{\varepsilon}}+\hat{B}(\mathbf{u^{\varepsilon}})-\mathbf{f}:=
-\nu\triangle\mathbf{u^{\varepsilon}}+[(\mathbf{u^{\varepsilon}}\cdot\nabla)
+(1/2)\Div\mathbf{u^{\varepsilon}}]\mathbf{u^{\varepsilon}}-\mathbf{f},$$
and $$\tilde{F}(\mathbf{u}):=A\mathbf{u} +
B(\mathbf{u})-\mathbf{f}:=-\nu\triangle\mathbf{u}+(\mathbf{u}\cdot\nabla)\mathbf{u}-\mathbf{f}.$$

 Let $\phi$ be a $\mathrm{C}^{\infty}$ scalar
function on $[0, T]$ with $\phi(T)=0$. Multiplying the equation
\eqref{e2.31} by $\phi(t)$, integrating in $t$ and taking
mathematical expectation, we obtain
\begin{align}\label{e2.35}
&(\mathbf{u^{\varepsilon}}(T),\mathbf{v})\phi(T) -
(\mathbf{u^{\varepsilon}}(0),\mathbf{v})\phi(0) -
E\Big[\int_{0}^{T}(\mathbf{u^{\varepsilon}}(t),\mathbf{v}\phi^{\prime}(t))\d
t\Big]\nonumber\\
&\quad + E\Big[\int_{0}^{T}
  \big(F(\mathbf{u^{\varepsilon}}(t)), \mathbf{v}\phi(t)\big)\d
  t\Big] + E\Big[\int_{0}^{T}\langle\nabla
p^{\varepsilon}(t),\mathbf{v}\phi(t)\rangle\, \d t\Big]\nonumber\\
&\qquad\quad = E\Big[\int_{0}^{T}(\mathbf{g}(t),\mathbf{v}\phi(t))\,
\d w(t)\Big], \quad\text{for all}\ \mathbf{v}\ \text{in}\
\mathbb{H}^1_0(\mathcal{O}).
\end{align}
Now passing to the limit in \eqref{e2.34} we have in the sense of
distribution,
\begin{align}
E\big[\sqrt{\varepsilon}\big(\frac{\d p^{\varepsilon}}{\d t},
q\big)\big]\longrightarrow E\big[\big(\frac{\d \chi}{\d t},
q\big)\big].
\end{align}
Hence in the same sense
\begin{align}
E\big[\varepsilon\big(\frac{\d p^{\varepsilon}}{\d t},
q\big)\big]\longrightarrow 0.\nonumber
\end{align}
Then passing to the limit in the following equation
\begin{align}
E\big[\langle \varepsilon \dot{p^{\varepsilon}}(t)+
\Div\mathbf{u^{\varepsilon}}(t),q\rangle\big] = 0\nonumber
\end{align}
we get
\begin{align}
E\big[\langle\Div\mathbf{u}, q\rangle\big] = 0,\quad \forall q \in
\mathrm{L}^{2}(\Omega;\mathrm{L}^{2}(0,T;\mathbb{L}^{2}(\mathcal{O}))),
\nonumber
\end{align}
which implies that $\Div\mathbf{u} = 0, \quad\text{almost everywhere
and almost
surely}$.\\
Hence $\mathbf{u}\in \mathrm{L}^{\mathrm{p}}(\Omega ;
  \mathrm{L}^{\infty}(0,T;\mathbb{L}^2(\mathcal{O}))) \cap \mathrm{L}^{2}
  (\Omega;\mathrm{L}^{2}(0,T;\mathbb{H}^1_0(\mathcal{O}))).$\\
Now it is clear that
\begin{align}\label{e2.37}
\langle\nabla p^{\varepsilon}, \mathbf{v}\phi(t)\rangle = \langle
p^{\varepsilon}, \phi(t)\Div\mathbf{v}\rangle = 0, \quad\text{almost
surely}.
\end{align}
Now using the same Minty-Browder monotonicity argument used in
Proposition 3.4 we can show that in limit
\begin{align}\label{e2.38}
F(\mathbf{u^{\varepsilon}})\longrightarrow
\tilde{F}(\mathbf{u})\quad\text{weakly in}\
\mathrm{L}^{2}(\Omega;\mathrm{L}^{2}(0,T;\mathbb{H}^{-1}(\mathcal{O}))).
\end{align}
Using the results \eqref{e2.37} and \eqref{e2.38}, in limit we have
from equation \eqref{e2.35}
\begin{align}
-E\Big[\int_{0}^{T}(\mathbf{u}(t),\mathbf{v}\phi^{\prime}(t))\d
t\Big] + E\Big[\int_{0}^{T}
  \big(\tilde{F}(\mathbf{u}(t)), \mathbf{v}\phi(t)\big)\d
  t\Big]= &(\mathbf{u}_0,\mathbf{v})\phi(0),\nonumber\\
& \text{for all}\ \mathbf{v}\ \text{in}\
\mathbb{H}^1_0(\mathcal{O}).\nonumber
\end{align}

This proves that $\mathbf{u}$ is a solution of the incompressible
stochastic Navier-Stokes equation.
\end{proof}

\begin{remark}
If one considers the multiplicative noise $\sigma(t,\mathbf{u})$ of
the type considered in the hypotheses $(A.1 - A.3)$ in Sritharan and
Sundar~\cite{SrSu06}, then under these conditions same a priori
estimates \eqref{e2.8}-\eqref{e2.10} hold. Thus
$\sigma^{\varepsilon}(.,\mathbf{u^{\varepsilon}})\rightarrow S$
weakly in $\mathrm{L}^{2}
  (\Omega;\mathrm{L}^{2}(0,T;\mathrm{L}_Q))$, where $\mathrm{L}_Q$
  denote the space of linear operators $S$ such that $SQ^{1/2}$ is a
  Hilbert-Schmidt operator from $\mathrm{L}^2$ to $\mathrm{L}^2$ and
  the norm on the space $\mathrm{L}_Q$ is defined by
  $|S|^{2}_{\mathrm{L}_Q}=\Tr (SQS^{\star})$ and $Q$ is a trace
  class operator. Hence with the help of Minty-Browder monotonicity arguments
  the existence and uniqueness proofs go through and we can also establish the limit
to incompressible flow.
\end{remark}

\par\bigskip\noindent
{\bf Acknowledgment.} We thank the reviewer for the helpful
comments.

\bibliographystyle{amsplain}

\end{document}